\documentclass[10pt]{article}
\usepackage{charter}\usepackage[T1]{fontenc}\usepackage{textcomp}
\usepackage[scaled=.92]{helvet}
\usepackage{eulervm,euscript}
\DeclareFontFamily{T1}{pzc}{}
\DeclareFontShape{T1}{pzc}{m}{it}{<-> s * [1.200] pzcmi}{}
\DeclareMathAlphabet{\nicemathcal}{OT1}{pzc}{m}{it}
\DeclareMathAlphabet{\oldmathcal}{OMS}{bch}{m}{n}
\usepackage{mathrsfs}
\usepackage{amsmath,amssymb,geometry,theorem}
\usepackage{xspace,calc,multirow,booktabs}
\usepackage[dvips]{graphicx}
\usepackage{pstricks,pst-node,pst-text,pst-tree,pstricks-add}
\psset{radius=2pt,xunit=1,yunit=1,arrowscale=1.7,arrowlength=.7,labelsep=3pt,rowsep=2\baselineskip,nodesep=.5ex}

\makeatletter
\def\section{\@startsection{section}{1}{0pt}{-3.25ex plus -1ex minus 
-.2ex}{1.5ex plus .2ex minus .3ex}{\normalfont\large\bf}}
\renewcommand\subsection{\@startsection{subsection}{2}{\z@}%
                                     {-3.25ex\@plus -1ex \@minus -.2ex}%
                                     {1.5ex \@plus .2ex}%
                                     {\normalfont\normalsize\bfseries}}
\makeatother
\renewenvironment{abstract}
{\vspace*{-1.8ex}\begin{quotation}\small
}{\end{quotation}}

\newcommand{\defn}[1]{{\textit{\textbf{#1}}}}
\newcommand{\myitem}[1]{\item[\textnormal{(#1)}]}

\theoremheaderfont{\scshape}
\theorembodyfont{\normalfont\slshape}
\theoremstyle{plain}

\newtheorem{lemma}{Lemma}

\newtheorem{corollary}{Corollary}
\newtheorem{theorem}{Theorem}
\theorembodyfont{\normalfont}

\newenvironment{proof}{\begin{trivlist}\item{}\normalfont\textit{Proof.}}{\hfill$\square$\end{trivlist}}
\newenvironment{proofof}[1]{\begin{trivlist}\item{}\normalfont\textit{Proof of #1.}}{\hfill$\square$\end{trivlist}}

\newcommand{\Ie}{\emph{I.e.}}
\newcommand{\ie}{\emph{i.e.}}
\newcommand{\eg}{\emph{e.g.}}
\newcommand{\qed}{\,\raisebox{0ex}{\rule{.8ex}{1.2ex}}}

\newdimen\arrayruleHwidth
\makeatletter
\def\Hline{\noalign{\ifnum0=`}\fi\hrule \@height \arrayruleHwidth
  \futurelet \@tempa\@xhline}
\makeatother

\newcommand{\uloop}[2]{\nccurve[ncurv=1.3,angleA=75,angleB=105]{#1}{#2}}
\newcommand{\dloop}[2]{\nccurve[ncurv=1.3,angleA=-75,angleB=-105]{#1}{#2}}
\newcommand{\uloopp}[2]{\nccurve[ncurv=1,angleA=75,angleB=105]{#1}{#2}}
\newcommand{\dloopp}[2]{\nccurve[ncurv=1,angleA=-75,angleB=-105]{#1}{#2}}

\newcommand{\lineangles}[4]{\nccurve[angleA=#3,angleB=#4]{#1}{#2}}

\newcommand{\vxrad}{1.7ex}
\newcommand{\vxgap}{\rule{\vxrad}{0ex}\rule{\vxrad}{0ex}}
\newcommand{\vx}[1]{\cnode*{2pt}{#1}}
\newcommand{\vxx}[1]{\cnode*{2pt}{#1}\vxgap}
\newcommand{\blankvx}[1]{\rnode{#1}{}}
\newcommand{\blankvxx}[1]{\rnode{#1}{}\vxgap}
\newcommand{\link}[1]{\cnode{2pt}{#1}}

\newcommand{\catfont}[1]{\textnormal{\textsf{#1}}}
\newcommand{\catopfont}[1]{\textnormal{\textbf{#1}}}

\renewcommand{\int}{\catopfont{Int}}

\newcommand{\Span}{\catopfont{Span}}

\newcommand{\iRel}{\catfont{iRel}}
\newcommand{\pfunc}{\catfont{pFun}}
\newcommand{\op}{^{\catopfont{op}}}
\newcommand{\Brau}{\catfont{Brau}}
\newcommand{\Part}{\catfont{Part}}

\newcommand{\TLieb}{\catfont{TLieb}}
\newcommand{\Link}{\catfont{Link}}
\newcommand{\flatLink}{\Link^{\flat}}
\newcommand{\flatPart}{\Part^{\flat}}

\newcommand{\flatBrau}{\Brau^{\flat}}
\newcommand{\flatTLieb}{\TLieb^{\flat}}
\newcommand{\MLL}{\catfont{MLL}}
\newcommand{\N}{\mathbb{N}}
\newcommand{\loopnode}{\cnode{\vxrad}{o}}
\newcommand{\powerset}{\mathcal{P}}
\newcommand{\linkinglen}[6]{\raisebox{-.3ex}{$\psset{labelsep=1pt}\Rnode X {#1}\hspace{#6} \Rnode W {#3} \hspace{#6} \Rnode Y {#5}\psset{arrows=->}
\ncline{W}{X}\taput{\,#2}\ncline{W}{Y}\taput{#4\,}$}}
\newcommand{\inlinelinkinglen}[6]{\raisebox{-.3ex}{$\psset{labelsep=1pt}\Rnode X {#1}\hspace{#6} \Rnode W {#3} \hspace{#6} \Rnode Y {#5}\psset{arrows=->}
\ncline{W}{X}\taput{\mbox{\scriptsize\,$#2$}}\ncline{W}{Y}\taput{\mbox{\scriptsize$#4$\,}}$}}
\newcommand{\inlinelinking}[5]{\inlinelinkinglen{#1}{#2}{#3}{#4}{#5}{4ex}}
\newcommand{\inlinecap}{\:\rnode{1}{\strut}\;\;\rnode{2}{\strut}\psset{nodesep=-.2pt}\uloop{1}{2}\:}
\newcommand{\inlinecup}{\:\rnode{1}{\strut}\;\;\rnode{2}{\strut}\psset{nodesep=-.2pt}\dloop{1}{2}\:}
\newcommand{\chmark}{{\hspace*{-.3ex}\mbox{\large\bf\checkmark}\hspace*{-.3ex}}}
\newcommand{\id}{\mathsf{i}}
\newcommand{\inrel}[3]{#1\,#3\,#2}
\newcommand{\inrels}[5]{#1\,#4\,#2\,#5\,#3}

\title{\vspace*{-3.5ex}\Large
\textbf{Linking diagrams for free}
\author{
   \large Dominic J.\ D.\ Hughes\thanks{Visiting Scholar, Computer Science Department, Stanford University, CA 94305, USA.}
   \\[1ex]
   \normalsize Stanford University
}
}
\date{\normalsize 9 May, 2008}

\begin{document}\thispagestyle{empty}
\maketitle

\begin{abstract}
Linking diagrams with path composition are ubiquitous, for example:
Temperley-Lieb and Brauer monoids, Kelly-Laplaza graphs for compact
closed categories, and Girard's multiplicative proof nets.  We
construct the category $\Link=\Span(\iRel)$, where $\iRel$ is the
category of injective relations (reversed partial functions) and show
that the aforementioned linkings, as well as Jones-Martin partition
monoids, reside inside $\Link$.  Path composition, including
collection of loops, is by pullback.  $\Link$ contains the free
compact closed category on a self-dual object (hence also the looped
Brauer and Temperly-Lieb monoids), and generalises partition monoids
with partiality (vertices in no partition) and empty- and infinite
partitions.  Thus we obtain conventional linking/partition diagrams
and their composition ``for free'', from $\iRel$.
\end{abstract}

\section{Introduction}\label{sec-intro}

Write $\flatBrau$ for the category of loopless Brauer linkings
\cite{Bra37}:
\begin{itemize}
\item \emph{Objects} $X,Y,\ldots$ are finite sets, 
whose elements we call \defn{vertices}.
\item \emph{Morphisms.} A \defn{loopless Brauer linking} $X\to Y$ 
is an equivalence relation on the disjoint union $X+Y$ whose every
class is a pair (2 vertices).
\item \emph{Composition} is path composition: the composite $SR:X\to Z$ of $R:X\to Y$ and $S:Y\to Z$ is the 
restriction to $X+Z$ of the transitive closure $(R+S)^*$ of
$R+S\subseteq X+Y+Z$.\footnote{\label{clutter}To avoid clutter we
assume here (without loss of generality, by renaming vertices) that
canonical injections $Q_i\to Q_1+Q_2$ are inclusions.  In other words,
we assume $X$, $Y$ and $Z$ are disjoint, and that every $+$ is a union
$\cup$.}
See Figure~\ref{fig-flatbrau-comp}.
\begin{figure*}
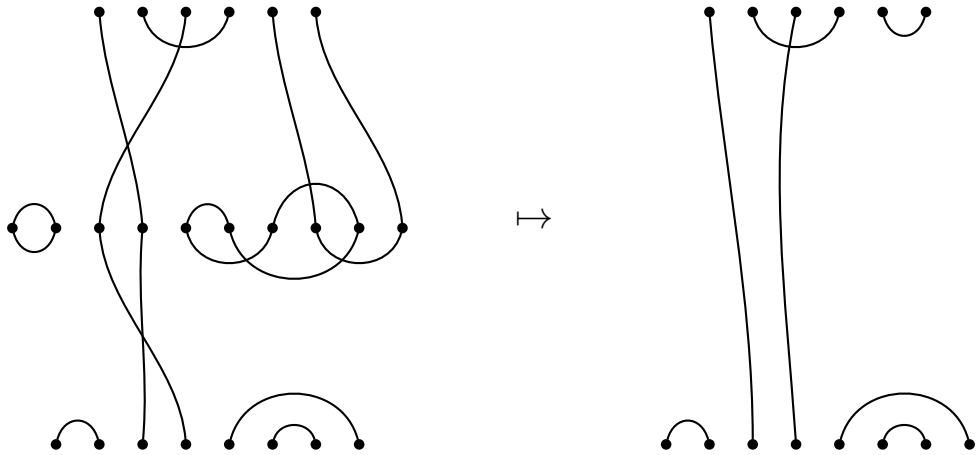

\begin{center}
\psset{nodesep=-.2mm}
\begin{math}
\begin{psmatrix}[rowsep=6ex]
\vxx{1}\vxx{2}\vxx{3}\vxx{4}\vxx{5}\vx{6} 
&&
\vxx{1a}\vxx{2a}\vxx{3a}\vxx{4a}\vxx{5a}\vx{6a} 
\\ 
\\
\vxx{1'}\vxx{2'}\vxx{3'}\vxx{4'}\vxx{5'}\vxx{6'}\vxx{7'}\vxx{8'}\vxx{9'}\vx{10'}
&\mbox{\Large$\mapsto$}&
\blankvxx{2'a}\blankvxx{4'a}\blankvxx{3'a}\blankvxx{5'a}\blankvxx{7'a}\blankvxx{8'a}\blankvxx{9'a}\blankvx{10'a}
\\
&&
\\
\vxx{1''}\vxx{2''}\vxx{3''}\vxx{4''}\vxx{5''}\vxx{6''}\vxx{7''}\vx{8''}
&&
\vxx{1''a}\vxx{2''a}\vxx{3''a}\vxx{4''a}\vxx{5''a}\vxx{6''a}\vxx{7''a}\vx{8''a}
\end{psmatrix}
\lineangles{1}{4'}{-85}{95}
\dloopp{2}{4}
\lineangles{3}{3'}{-95}{85}
\lineangles{5}{8'}{-85}{95}
\lineangles{6}{10'}{-85}{95}
\uloop{1'}{2'}
\uloop{5'}{6'}
\uloop{7'}{9'}
\dloop{1'}{2'}
\dloopp{5'}{7'}
\dloopp{6'}{9'}
\dloopp{8'}{10'}
\lineangles{3'}{4''}{-85}{95}
\lineangles{4'}{3''}{-95}{85}
\uloop{1''}{2''}
\uloopp{5''}{8''}
\uloopp{6''}{7''}
\lineangles{1a}{3''a}{-85}{90}
\dloopp{2a}{4a}
\lineangles{3a}{4''a}{-102}{94}
\dloop{5a}{6a}
\uloop{1''a}{2''a}
\uloopp{5''a}{8''a}
\uloopp{6''a}{7''a}
\end{math}
\end{center}
\caption{\label{fig-flatbrau-comp}Example of composition in 
the category $\flatBrau$ of loopless Brauer linkings.  Each
equivalence class $\{x,y\}$ is depicted as a ``link'' on $x$ and $y$.}
\end{figure*}
\end{itemize}
The \defn{loopless Brauer monoid} $\flatBrau_n$ is the subcategory of
$\flatBrau$ on $\{1,\ldots,n\}$.\footnote{\Ie, the monoid
$\flatBrau_n$ is the homset
$\flatBrau(\{1,\ldots,n\},\{1,\ldots,n\})$, with composition as
multiplication.  Although \cite{Bra37} considered only monoids,
collecting them into a category is obvious and trivial.}

Write $\Brau$ for the category of looped Brauer linkings, on the same objects:
\begin{itemize}
\item \emph{Morphisms.} A \defn{looped Brauer linking}
$X\to Y$ is a pair $\langle k,R\rangle$, denoted $\delta^k R$,
comprising a loopless Brauer linking $R:X\to Y$ a \defn{loop count}
$k\in\N=\{0,1,\ldots\}$.
\item \emph{Composition} is path composition, collecting loops: $(\delta^l S)\,(\delta^k R)$ 
is $\delta^{l+k+\lambda}(SR)$ where $SR$ is the composite in
$\flatBrau$ and $\lambda$ is the number of \defn{loops} formed during
the construction of $SR$, that is, classes of $(R+S)^*\subseteq X+Y+Z$
which are entirely within $Y$.
See Figure~\ref{fig-composition}.
\begin{figure*}
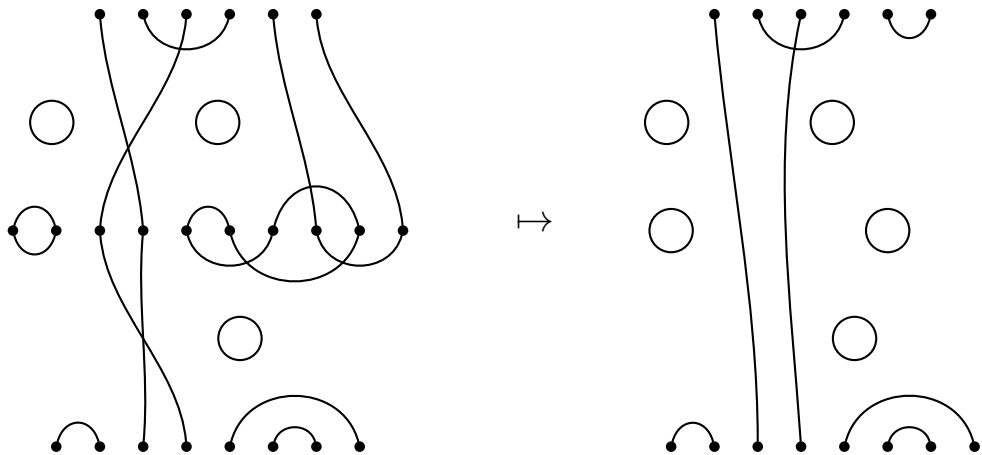

\begin{center}
\psset{nodesep=-.2mm}
\begin{math}
\begin{psmatrix}[rowsep=6ex]
\vxx{1}\vxx{2}\vxx{3}\vxx{4}\vxx{5}\vx{6} 
&&
\vxx{1a}\vxx{2a}\vxx{3a}\vxx{4a}\vxx{5a}\vx{6a} 
\\ 
\loopnode\hspace{13ex}\loopnode\hspace{11.5ex}\mbox{}
&&
\loopnode\hspace{13ex}\loopnode\hspace{11.5ex}\mbox{}
\\
\vxx{1'}\vxx{2'}\vxx{3'}\vxx{4'}\vxx{5'}\vxx{6'}\vxx{7'}\vxx{8'}\vxx{9'}\vx{10'}
&\mbox{\Large$\mapsto$}&
\loopnode\blankvxx{2'a}\blankvxx{4'a}\blankvxx{3'a}\blankvxx{5'a}\blankvxx{7'a}\loopnode\blankvxx{8'a}\blankvxx{9'a}\blankvx{10'a}
\\
\loopnode\hspace{-5ex}\mbox{}
&&
\loopnode\hspace{-5ex}\mbox{}
\\
\vxx{1''}\vxx{2''}\vxx{3''}\vxx{4''}\vxx{5''}\vxx{6''}\vxx{7''}\vx{8''}
&&
\vxx{1''a}\vxx{2''a}\vxx{3''a}\vxx{4''a}\vxx{5''a}\vxx{6''a}\vxx{7''a}\vx{8''a}
\end{psmatrix}
\lineangles{1}{4'}{-85}{95}
\dloopp{2}{4}
\lineangles{3}{3'}{-95}{85}
\lineangles{5}{8'}{-85}{95}
\lineangles{6}{10'}{-85}{95}
\uloop{1'}{2'}
\uloop{5'}{6'}
\uloop{7'}{9'}
\dloop{1'}{2'}
\dloopp{5'}{7'}
\dloopp{6'}{9'}
\dloopp{8'}{10'}
\lineangles{3'}{4''}{-85}{95}
\lineangles{4'}{3''}{-95}{85}
\uloop{1''}{2''}
\uloopp{5''}{8''}
\uloopp{6''}{7''}
\lineangles{1a}{3''a}{-85}{90}
\dloopp{2a}{4a}
\lineangles{3a}{4''a}{-102}{94}
\dloop{5a}{6a}
\uloop{1''a}{2''a}
\uloopp{5''a}{8''a}
\uloopp{6''a}{7''a}
\end{math}
\end{center}
\caption{\label{fig-composition}Example
of composition in the category $\Brau$ of looped Brauer linkings.  The
two input linkings are $\delta^2 R$ (upper) and $\delta^1 S$ (lower),
where $R$ and $S$ are the loopless linkings in
Figure~\ref{fig-flatbrau-comp}.  The output linking $(\delta^1
S)\,(\delta^2 R)$ is $\delta^5 (SR)\,=\,\delta^{1+2+2}\:(SR)$, where
$SR$ is the output loopless linking in Figure~\ref{fig-flatbrau-comp},
a composition which forms two new loops.}
\end{figure*}
\end{itemize}
The \defn{looped Brauer monoid} $\Brau_n$ is the subcategory of
$\Brau$ on $\{1,\ldots,n\}$.\footnote{$\Brau_n$ is the submonoid of
the Brauer algebra over $n$ \cite{Bra37} generated (under
multiplication in the algebra) by $\{\delta^0 R\,:\,R\in\flatBrau_n\}$
and $\delta^1\id$, where $\id$ is the identity in $\flatBrau_n$.}
The category $\Brau$ is (equivalent to) the free compact closed
category on a self-dual object \cite{KL80,Abr05}.
There is a forgetful functor to both $\Brau$ and $\flatBrau$ from the
category $\MLL$ of unit-free multiplicative proof nets \cite{Gir87},
extracting leaves (literal occurences) and axiom links.\footnote{An
object of $\MLL$ is a unit-free multiplicative formula, a morphism
$A\to B$ is a cut-free proof net on $A\multimap B$, and composition is by
cut elimination.  See \eg\ \cite{HG03,HG05}.
The well-definedness to $\flatBrau$ is trivial; the functor to $\Brau$
is more subtle, being well-defined because proof net correctness
ensures no loops arise during composition (\ie, $\lambda=0$ in the
definition of composition in $\Brau$).}

The separate treatment of paths and loops is ad hoc.  We shall unify
paths and loops, handling them simultaneously, and in so doing, obtain
infinite generalisations of linkings.

\paragraph*{Acknowledgement.}  Thanks to Robin Houston for feedback
last summer on the prospect of extending pullbacks from injective
relations to coherence spaces \cite{Gir87} for a ``sliced'' notion of
linking, enriched in commutative monoids.
This is work in progress.

Many thanks to Vaughan Pratt for his ongoing support.

\section{Generalised linkings: $\Link=\Span(\iRel)$}

A binary relation $R:A\to Z$ (\ie, $R\subseteq A\times Z$) is
\defn{injective} if $aRz$ and $a'Rz$ implies $a=a'$.\footnote{$aRz$
abbreviates $\langle a,z\rangle\in R$.}
Write $\iRel$ for the category of sets and injective relations between
them.  Note that $\iRel=\pfunc\op$, the opposite of the category of
sets and partial functions.

\begin{figure*}
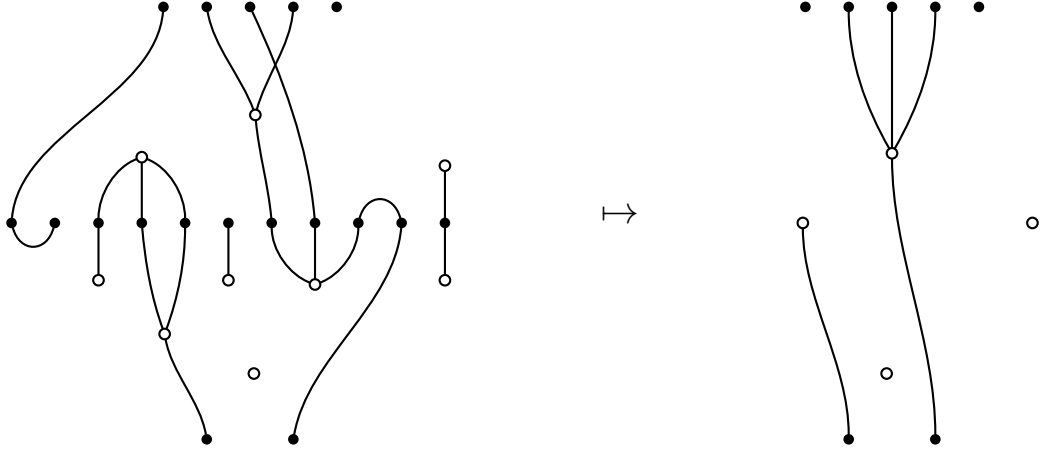

\begin{center}
\psset{nodesep=-.2mm,labelsep=3pt}
\begin{math}
\begin{psmatrix}[rowsep=6ex]
\vxx{1}\vxx{2}\vxx{3}\vxx{4}\vx{5} 
&&
\vxx{1a}\vxx{2a}\vxx{3a}\vxx{4a}\vx{5a} 
\\ 
\\
\vxx{1'}\vxx{2'}\vxx{3'}\vxx{4'}\vxx{5'}\vxx{6'}\vxx{7'}\vxx{8'}\vxx{9'}\vxx{10'}\vxx{11'}
&\mbox{\Large$\mapsto$}&
\hspace*{4ex}\link{ww}\hspace*{18ex}\link{w}
\\
\\
\vxx{1''}\vxgap\vx{2''}
&&
\vxx{1''a}\vxgap\vx{2''a}
\end{psmatrix}
\lineangles{1}{1'}{-95}{85}
\dloop{1'}{2'}
\nccurve[linestyle=none]{2}{5'}\ncput[npos=.5]{\link{w1}}
\lineangles{2}{w1}{-80}{110}
\lineangles{w1}{7'}{-85}{95}
\lineangles{4}{w1}{-95}{75}
\lineangles{3}{8'}{-65}{95}
\ncline[linestyle=none]{3'}{5'}\naput[labelsep=5ex]{\link{w2}}
\lineangles{w2}{3'}{-160}{90}
\lineangles{w2}{4'}{-90}{90}
\lineangles{w2}{5'}{-20}{90}
\nput[labelsep=4ex]{-90}{3'}{\link{w3}}
\ncline{w3}{3'}
\nput[labelsep=4ex]{-90}{6'}{\link{w4}}
\ncline{w4}{6'}
\ncline[linestyle=none]{4'}{1''}\nbput{\link{w5}}
\lineangles{w5}{4'}{110}{-85}
\lineangles{w5}{5'}{70}{-90}
\lineangles{w5}{1''}{-80}{100}
\ncline[linestyle=none]{7'}{9'}\naput[labelsep=-5ex]{\link{w6}}
\lineangles{w6}{7'}{160}{-90}
\lineangles{w6}{8'}{90}{-90}
\lineangles{w6}{9'}{20}{-90}
\uloop{9'}{10'}
\lineangles{10'}{2''}{-95}{80}
\nput[labelsep=4ex]{90}{11'}{\link{w7}}
\ncline{w7}{11'}
\nput[labelsep=4ex]{-90}{11'}{\link{w8}}
\ncline{w8}{11'}
\nput[labelsep=11ex]{-90}{3a}{\link{w}}
\lineangles{w}{2a}{120}{-90}
\lineangles{w}{3a}{90}{-90}
\lineangles{w}{4a}{60}{-90}
\lineangles{w}{2''a}{-90}{90}
\lineangles{ww}{1''a}{-90}{90}
\ncline[linestyle=none]{1''}{2''}\naput[npos=.55,labelsep=5ex]{\link{.}}
\ncline[linestyle=none]{1''a}{2''a}\naput[npos=.43,labelsep=5ex]{\link{.}}
\end{math}
\end{center}
\caption{\label{fig-linkings}Examples of linkings and pullback-composition in $\Link=\Span(\iRel)$.  A link is shown 
as a small circle, with its vertices attached by edges.
We leave the circle implicit when a link has two vertices.}
\vspace{1ex}\hrule
\end{figure*}

A \defn{linking} $X\to Y$ is a diagram
\begin{equation*}
\rule{0ex}{3ex}
\Rnode X X\hspace{5ex} \Rnode A A \hspace{5ex} \Rnode Y Y\psset{arrows=->}
\ncline{A}{X}\taput{f}\ncline{A}{Y}\taput{g}
\end{equation*}
in $\iRel$.  Each $a\in A$ is a \defn{link},
and the elements of $X$ and $Y$ are \defn{vertices}.
The vertex set $f(a)+g(a)\subseteq X+Y$ is the \defn{footprint} of the link
$a\in A$.\footnote{For any binary relation $R:A\to Z$, the image
$R(a)$ is $\{z\in Z:aRz\text{ for some }a\in A\}\subseteq Z$.}
If a vertex $x$ is in the footprint of $a$, we simply say that $x$ is
in $a$, or $a$ has/contains $x$.
The injectivity requirement implies that no two links overlap (share a vertex).
See Figure~\ref{fig-linkings} for examples.

Just as graph theory treats graphs up to isomorphism, we identify
linkings up to isomorphism, \ie, renaming of links.
Formally, we identify linkings $\inlinelinking{X}{f}{A}{g}{Y}$ and
$\inlinelinking{X}{\,f'}{A'}{\;\;\;g'}{Y}$ iff there exists a
bijection $\theta:A\to A'$ such that $f'\theta=f$ and
$g'\theta=g$.
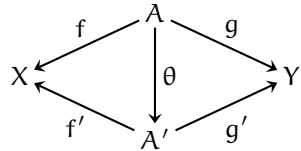
$$\psset{arrows=->,labelsep=2pt}
\begin{array}{c}
\Rnode A A \\[2.5ex]
\Rnode X X\hspace{20ex} \Rnode Y Y \\[2.5ex]
\Rnode{A'}{A'}
\end{array}
\ncline{A}{X}\nbput{f}\ncline{A}{Y}\naput{g}
\ncline{A'}{X}\naput{f'}\ncline{A'}{Y}\nbput{g'}
\ncline{A}{A'}\trput{\theta}
$$

\subsection{Composition by pullback}
The composite $X\to Z$ of linkings $\inlinelinking{X}{f}{A}{g}{Y}$ and
$\inlinelinking{Y}{\,h}{B}{k}{Z}$ is by pullback in
$\iRel$:\footnote{Equivalently, pushout in $\pfunc=\iRel\op$.  The
use of spans/pullbacks in this paper, together with equivalence up to
isomorphism, should compared with the standard use of cospans/pushouts
for tangles and cobordisms.}
\begin{equation}\label{diag-pullback}
\newcommand{\hozgap}{\hspace{13ex}}\psset{arrows=->,labelsep=2pt}
\begin{array}{c}
\Rnode P P \\[4ex]
\Rnode A A \hozgap \Rnode{B}{B} \\[4ex]
\Rnode X X \hozgap \Rnode Y Y \hozgap \Rnode Z Z
\end{array}
\ncline{A}{X}\nbput{f}\ncline{A}{Y}\naput{g}
\ncline{B}{Y}\nbput{h}\ncline{B}{Z}\naput{k}
\ncline{P}{A}\nbput{p}
\ncline{P}{B}\naput{q} 
\end{equation}
Explicitly, the composite linking
$$\linkinglen{X}{fp}{P}{k\mkern-2mu q}{Z}{8ex}$$
is defined as follows. 
To illustrate the definition as we proceed, we refer to the Brauer
composition in Figure~\ref{fig-composition}.  There $X/Y/Z$ are the
upper/mid/lower rows, and $A/B$ are the upper/lower link sets.

A \defn{synchronisation} $\langle \alpha,\beta\rangle$ is a pair of
sets of links $\alpha\subseteq A$ and $\beta\subseteq B$ with the same
footprint in the interface $Y$:
\begin{equation}\label{eq-sync}
g(\alpha)\;\;\;=\;\;\;h(\beta)
\end{equation}
For example, in Figure~\ref{fig-composition}, if $\alpha$ comprises
the three caps $\inlinecap$ of $A$, and $\beta$ the first three cups
$\inlinecup$ of $B$, then $\langle
\alpha,\beta\rangle$ is a synchronisation with
$f(\alpha)=g(\alpha)=\{y_1,y_2,y_5,y_6,y_7,y_9\}\subseteq Y$, where
the $y_i$ are the vertices of $Y$ from left to right.  (Note that this
remains a synchronisation upon adding any number of loops to $\alpha$
and $\beta$, since loops have empty footprint in $Y$.)  
Henceforth identify a synchronisation $\langle \alpha,\beta \rangle$
(and more generally any pair $\langle \alpha,\beta\rangle$ of subsets
$\alpha\subseteq A$ and $\beta\subseteq B$) with the corresponding
subset $\alpha+\beta\subseteq A+B$ (thus identifying along the
bijection\footnote{More suggestively, $2^A\times 2^B\,\cong\,2^{A+B}$,
writing $2^C$ for $\powerset(C)$.}
$\powerset(A)\times\powerset(B)\,\cong\,\powerset(A+B)$, where
$\powerset(C)$ denotes the powerset (set of subsets) of $C$).

A (generalised) \defn{path} is a minimal non-empty synchronisation,
where minimality is with respect to inclusion.
There are 12 paths in Figure~\ref{fig-composition}: seven singletons
(the two loops in $A$, the loop in $B$, the cup of $A$, and the three
caps of $B$), three doubletons (the short circuit formed on
$\{y_1,y_2\}$ and the verticals through $y_3$ and $y_4$), one triplet
(through $y_8$ and $y_{10}$), and one quadruplet (the long circuit through
$y_5,y_6,y_7,y_9$).

Define the set $P$ of links of the composite 
$\inlinelinkinglen{X}{fp}{P}{k\mkern-2mu q}{Z}{5ex}$
as the set of all paths, and
define $p:P\to A$ and $q:P\to B$ as the projections
\begin{align}
p\langle\alpha,\beta\rangle &= \alpha\label{eq-p} \\
q\langle\alpha,\beta\rangle &= \beta\label{eq-q}
\end{align}
In Figure~\ref{fig-composition}, $p$ (resp.\ $q$) projects each path
to its constituent links in the upper half $A$ (resp.\ lower half $B$).
The composite $fp:P\to X$ projects a path $\gamma$ to the vertices (if
any) in $X$ which are on $\gamma$, and similarly for $kq:P\to Y$.
In particular, for each of the five loops $L$ (both the three
singletons from the original linkings, and the two formed of multiple
links), we have $fp(L)$ and $kq(L)$ empty.

See Figure~\ref{fig-linkings} for a more general, non-Brauer example.
An example of an infinite composition is depicted in
Figure~\ref{fig-inf-comp}, illustrating why naive infinite
generalisations of Brauer linkings do not work: an infinite chain of
binary (two-vertex) links produces a unary (single-vertex) link.
A finite variant is in Figure~\ref{fig-fin-comp}.
\begin{figure*}
\begin{center}
\psset{nodesep=-.2mm,labelsep=3pt}
\begin{math}
\begin{psmatrix}[rowsep=4ex]
\vxx{1}\vxgap\vxgap\vxgap\vxx{3}\vxx{4}\vxgap 
&&
\vxx{1a}\vxgap\vxgap\vxgap\vxx{3a}\vxx{4a}\vxgap
\\ 
\\
\vxx{1'}\vxx{2'}\vxx{3'}\vxx{4'}\vxx{5'}\vxx{6'}\vxx{7'}\blankvx{8'}\makebox[0ex]{\raisebox{-.2ex}{$\;\;\;\;\;\ldots$}}
&\mbox{\Large$\mapsto$}&
\link{w}\vxgap\vxgap\vxgap\vxgap\vxgap\vxgap\vxgap
\\
\\
\vxx{1''}\vxx{2''}\vxx{3''}\vxx{4''}
&&
\vxx{1''a}\vxx{2''a}\vxx{3''a}\vxx{4''a}
\end{psmatrix}
\ncline{1}{1'}
\dloop{1'}{2'}
\uloop{2'}{3'}
\dloop{3'}{4'}
\uloop{4'}{5'}
\dloop{5'}{6'}
\uloop{6'}{7'}
{\psset{nodesepB=.5ex}\dloop{7'}{8'}}
\dloop{3}{4}
\dloop{3a}{4a}
\uloop{1''}{2''}
\uloop{3''}{4''}
\uloop{1''a}{2''a}
\uloop{3''a}{4''a}
\ncline{w}{1a}
\end{math}
\end{center}
\caption{\label{fig-inf-comp}Example
of the pullback-composition of infinite linkings in
$\Link=\Span(\iRel)$, with interface layer $Y=\{1,2,\ldots\}$.  The
entirety of $Y$ is a synchronisation, and since it is minimal and
non-empty, it is a path.  Thus it shows up in the result of
composition on the right, as a unary link (with a single vertex).
This shows clearly why naive infinite generalisations of Brauer
linkings do not work: an infinite chain of binary (two-vertex) links
has produced a unary link.
Figure~\ref{fig-fin-comp} shows a corresponding example in which the
interface layer is finite.}
\vspace*{8ex}
\begin{center}
\psset{nodesep=-.2mm,labelsep=3pt}
\begin{math}
\begin{psmatrix}[rowsep=4ex]
\vxx{1}\vxgap\vxgap\vxgap\vxx{3}\vxx{4}\vxgap 
&&
\vxx{1a}\vxgap\vxgap\vxgap\vxx{3a}\vxx{4a}\vxgap
\\ 
\\
\vxx{1'}\vxx{2'}\vxx{3'}\vxx{4'}\vxx{5'}\vxx{6'}\vx{7'}\vxgap
&\mbox{\Large$\mapsto$}&
\\
\\
\vxx{1''}\vxx{2''}\vxx{3''}\vxx{4''}
&&
\vxx{1''a}\vxx{2''a}\vxx{3''a}\vxx{4''a}
\end{psmatrix}
\ncline{1}{1'}
\dloop{1'}{2'}
\uloop{2'}{3'}
\dloop{3'}{4'}
\uloop{4'}{5'}
\dloop{5'}{6'}
\uloop{6'}{7'}
\dloop{3}{4}
\dloop{3a}{4a}
\uloop{1''}{2''}
\uloop{3''}{4''}
\uloop{1''a}{2''a}
\uloop{3''a}{4''a}
\end{math}
\end{center}
\caption{\label{fig-fin-comp}Analogous
composition to Figure~\ref{fig-inf-comp} in which the interface layer
$Y=\{1,\ldots,7\}$ is finite.  This time there is no non-empty
synchronisation touching $Y$, so no link results from the interaction
there.  (Note that the lower-left input linking is not a Brauer linking,
since it is partial: vertex 7 is in no link.)}
\end{figure*}

\begin{theorem}\label{thm-pullback}
The construction above defines pullbacks in\/ $\iRel$.
\end{theorem}
\begin{proof}
Section~\ref{sec-pullback}.
\end{proof}
Write $\Link$ for the category of linkings with this composition.  In
other words, $\Link=\Span(\iRel)$, the span construction \cite{Ben67}
applied to $\iRel$, with bicategorical structure collapsed to a
category by taking morphisms (1-cells) up to isomorphism.  That
$\Link$ is a category (with identities and associative composition)
follows from the general features of the $\Span$ construction, saving
considerable labour.

\subsection{Loopless variant $\flatLink$}

A \defn{loop} is a link without vertices.  Define $\flatLink$ as the
variant of $\Link$ comprising the loopless linkings, discarding any
loops formed during pullback composition.
(Composition is associative since loops do not interact during
pullback.)
Write $(-)^\flat :\Link\to\flatLink$ for the functor which deletes
loops (identity on objects).
Note that $\flatLink$ is not a subcategory of $\Link$, since composition of
loopless linkings can generate loops.

\section{Subcategories of $\Link$ and $\flatLink$}

We consider various subcategories of $\Link$ and $\flatLink$, as
summarised in Figure~\ref{fig-subcats} and detailed below.\begin{figure*}
\begin{center}
\begin{equation*}
\begin{array}{c@{}c@{\;\;\subseteq\;\;}c@{\;\;\subseteq\;\;}c@{\;\;\subseteq\;\;}c}
(\N,+) \;\;\subseteq\;\; & \Rnode{tlieb}{\TLieb} & \Rnode{brau}{\Brau} & \Rnode{part}{\Part} 
& \Rnode{link}{\Link} \\[5ex]
& \Rnode{flattlieb}{\flatTLieb} & \Rnode{flatbrau}{\flatBrau} & \Rnode{flatpart}{\flatPart} 
& \Rnode{flatlink}{\flatLink}
\end{array}
\psset{arrows=->,nodesepA=5pt}
\ncline{tlieb}{flattlieb}\trput{\flat}
\ncline{brau}{flatbrau}\trput{\flat}
\ncline{part}{flatpart}\trput{\flat}
\ncline{link}{flatlink}\trput{\flat}
\end{equation*}
\vspace{6ex}

\begin{tabular}{ccccccc}
\toprule
& \multirow{2}{*}{\begin{tabular}{c}Object\\
restriction\\[-1.5ex]\end{tabular}} & \multicolumn{5}{c}{Morphism
restriction} \\
\cmidrule(l){3-7}
& & \it loopless & \it finite & \it total & \it binary & \it planar \\
\midrule
\rule{0ex}{2.5ex}%
$\Link$\\[.5ex]
$\Part$  & finite & & \chmark & \chmark \\[.5ex]
$\Brau$  & finite & & \chmark & \chmark & \chmark \\[.5ex]
$\TLieb$ & \{1,\ldots,n\} & & \chmark & \chmark & \chmark & \chmark  
\\[.5ex]
$(\N,+)$     & empty & & \chmark & (\chmark) & (\chmark) & (\chmark)
\\[1ex]
\midrule
\rule{0ex}{2.5ex}%
$\flatLink$  &         & \chmark \\[.5ex]
$\flatPart$  & finite & \chmark & (\chmark) & \chmark \\[.5ex]
$\flatBrau$  & finite & \chmark & (\chmark) & \chmark & \chmark \\[.5ex]
$\flatTLieb$ & \{1,\ldots,n\} & \chmark & (\chmark) & \chmark & \chmark & \chmark  
\\[1ex]
\bottomrule
\end{tabular}
\vspace{3ex}
\end{center}
\caption{\label{fig-subcats}Various subcategories of $\Link$ and its 
loopless variant $\flatLink$.  \emph{Total} means every vertex is in a
link.  \emph{Binary} means every non-loop has exactly two vertices.
The $(\chmark)$ are implied $\chmark$.  Here $(\N,+)$ is the monoid of
integers under addition, which is the subcategory of $\Link$ on the
empty set.  The functor $(-)^\flat$ deletes all loops. The categories
$\flatPart$, $\flatBrau$ and $\flatTLieb$ contain the standard
(loopless) Jones-Martin partition-, Brauer- and Temperley-Lieb monoids,
respectively.}
\end{figure*}%

The categories $\Brau$ and $\flatBrau$ were defined at the start of
Section~\ref{sec-intro}.
The categories $\Part$ and $\flatPart$ are the looped and unlooped
\defn{Jones-Martin partition categories}
\cite{Jon94,Mar94}\footnote{\label{obvious}As with the Brauer
category, we have merely collected the monoids into categories in the
obvious way.}, defined exactly as $\Brau$ and $\flatBrau$ (verbatim),
but dropping the restriction that every equivalence class is a pair.
The conventional (loopless) partition monoid on $n$ is the subcategory
$\flatPart_n$ of $\flatPart$ on $\{1,\ldots,n\}$.

The \defn{Temperley-Lieb category} $\TLieb$ \cite{TL71}\footnote{See
footnote~\ref{obvious}.} is the subcategory of $\Brau$ on objects of
the form $\{1,\ldots,n\}$ for $n\ge 0$, and with only the
planar\footnote{We assume vertices $1,\ldots,n$ are ordered in the
plane.}  linkings (no crossings of links,
\ie, well-bracketed or ``parenthetical'' \cite[p.\,63]{Kau04}).  See
\cite{Abr07a} for a concrete presentation.
The category $\flatTLieb$ is the loopless variant of $\TLieb$.  The
standard loopless Temperley-Lieb monoids are the subcategories of
$\flatTLieb$ on the objects $\{1,\ldots,n\}$.

Planar partition monoids can be defined by analogy with Temperley-Lieb
monoids.  For a nice exposition of each of the aforementioned monoids
(and their algebras), with many diagrams and examples, see
\cite{HR05}.

\section{Geometry of interaction ``for free''}

Let $\MLL$ denote the category of multiplicative proof nets
\cite{Gir87}, with unit-free formulas as objects, a morphism $X\to Y$
as a cut-free proof net on $X\multimap Y$, and composition by cut
elimination.  Thus a proof net is a linking on leaves (literal
occurrences) which satisfies a correctness criterion, and composition
is path composition.\footnote{See \eg\ \cite{HG03,HG05}.}  The
forgetful functor $L\!^\flat:\MLL\to\flatBrau$ extracts the leaves
(forgetting the underlying parse tree structure of the formulas) and
the links between them.
Due to the correctness criterion on proof nets, loops never arise
during composition, thus there is also a forgetful functor
$L:\MLL\to\Brau$, and the following diagram commutes.
\begin{equation*}
\begin{array}{c@{\hspace{5ex}}c@{\hspace{1ex}\subseteq\hspace{1ex}}c}
\Rnode{mll}{\MLL} & \Rnode{brau}{\Brau} & \Rnode{link}{\Link} \\[6ex]
& \Rnode{fbrau}{\flatBrau} & \Rnode{flink}{\flatLink}
\end{array}
\psset{arrows=->,nodesepA=5pt}
\ncline{mll}{brau}\taput{L}
\ncline{mll}{fbrau}\nbput[npos=.4]{L\!^\flat}
\ncline{brau}{fbrau}\trput{\flat}
\ncline{link}{flink}\trput{\flat}
\end{equation*}
Having composed the linkings of proof nets $A\multimap B$ and $B\multimap C$
in $\Link$ by $\iRel$ pullback, we can draw the resulting linking on
$A\multimap C$, to obtain the composite in $\MLL$.  Thus all computation
happens inside $\Link$, so we have geometry of interaction
\cite{Gir89} ``for free'', via $\iRel$.

Work in progress aims to use pullbacks of coherence spaces
\cite{Gir87}, an extension of $\iRel$, to obtain a
multiplicative-additive geometry of interaction ``for free''.

\section{Proof of Theorem~\ref{thm-pullback}}\label{sec-pullback}

A binary relation $R:A\to Z$ is \defn{total} if the image
$R(a)\subseteq Z$ is non-empty for all $a\in A$.
\begin{lemma}\label{lemma-monic}
An $\iRel$ morphism is monic\footnote{Recall that a morphism $m:A\to Z$ is
\emph{monic} if $mf=mg$ implies $f=g$ for all objects $W$ and
$f,g:W\to A$ \cite{Mac71}.} iff it is total.\footnote{Dually, and
perhaps more intuitively obvious, a partial function is epic (in
$\pfunc$) iff it is surjective.}
\end{lemma}
\begin{proof}
Suppose $m:A\to Z$ is total.  Let $f,g:W\to A$ with $mf=mg$.  If
$f\neq g$ there exist $w\in W$ and $a\in A$ with $\inrel{w}{a}{f}$ but
not $\inrel w a g$ (exchanging $f$ and $g$, if necessary).  Since $m$
is total, there exists $z\in Z$ with $\inrel{a}{z}{m}$.  Thus
$\inrel{w}{z}{(mf)}$, so $\inrel{w}{z}{(mg)}$, hence there exists
$a'\in A$ with $\inrels{w}{a'}{z}{g}{m}$.  Since not
$\inrel{w}{a}{g}$, we have $a'\!\neq\! a$, but then $\inrel a z m$ and
$\inrel {a'} z m$ contradicting injectivity.  Thus $f\!=\!g$, so $m$ is
monic.

Conversely, suppose $m:A\to Z$ is not total.  Then there exists $a\in
A$ such that $m(a)=\emptyset$.  Let $W=\{w\}$, $\,f(w)=\emptyset$ and
$\,g(w)=\{a\}\,$.  Then $mf=mg$ (both empty) yet $f\!\neq\! g$, so $m$ is not monic.
\end{proof}

\begin{lemma}[Stability]
Injective relations preserve unions and intersections: for any $R:A\to
Z$ in\/ $\,\iRel$ and subsets $\alpha_i\subseteq A$ for each $i$ in some indexing set $I$,
\begin{align}
\label{eq-stab-cup}R\left(\:\bigcup_{i\in I}\alpha_i\right)\;\;\;=\;\;\;\bigcup_{i\in i}\,R(\alpha_i) \\[1ex]
\label{eq-stab-cap}R\left(\:\bigcap_{i\in I}\alpha_i\right)\;\;\;=\;\;\;\bigcap_{i\in i}\,R(\alpha_i)
\end{align}
\end{lemma}
\begin{proof}
\eqref{eq-stab-cup}. A trivial property of binary relations (injectivity not required).\footnote{$b\in R(\bigcup\alpha_i)$ 
iff $\exists a\in\bigcup\alpha_i.aRb$ iff $\exists i\in I,
a\in\alpha_i.aRz$ iff $\exists i\in I.z\in R(\alpha_i)$ iff $z\in
\bigcup R(\alpha_i)$.}

\eqref{eq-stab-cap}.  Suppose $z\in R(\bigcap\alpha_i)$, \ie, $aRz$ for some
$a\in\bigcap\alpha_i$.  Then $a\in\alpha_i$ for all $i$, hence $z\in
R(\alpha_i)$ for all $i$, so $z\in
\bigcap R(\alpha_i)$.
Conversely, suppose $z\in\bigcap R(\alpha_i)$,
\ie, $z\in R(\alpha_i)$ for all $i$.  Then for each $i\in I$ there
exists $a_i\in \alpha_i\subseteq A$ with $a_iR z$.  By injectivity,
$a_i=a_j=a$ for all $i,j\in I$, hence $a\in\bigcap \alpha_i$.  Thus
$z\in R(\bigcap\alpha_i)$, since $aRz$.
\end{proof}
Write $\alpha\uplus\beta$ for $\alpha\cup\beta$ when
$\alpha\cap\beta=\emptyset$, and more generally, write
$\biguplus_{i\in I}\alpha_i$ for $\bigcup_{i\in I}\alpha_i$ when
$\alpha_i\cap\alpha_j=\emptyset$ for all distinct $i,j\in I$.
\begin{corollary}\label{cor-stab-uplus}
Injective relations preserve disjoint unions: with $R$ as in the previous lemma,
\begin{align}
\label{eq-stab-uplus}R\left(\:\biguplus_{i\in I}\alpha_i\right)\;\;\;=\;\;\;\biguplus_{i\in i}\,R(\alpha_i)
\end{align}
\end{corollary}
\begin{proof}
Immediate from \eqref{eq-stab-cup} and \eqref{eq-stab-cap}.
\end{proof}
\begin{corollary}\label{cor-irel-sub}
Injective relations preserve
inclusion and
subtraction:\/
if\/ $R:A\to Z$ in\/ $\iRel$ and $\alpha,\beta\subseteq A$ then\footnote{$\alpha\setminus\beta\;=\;\{a\in\alpha:a\not\in\beta\}$.}
\begin{eqnarray}
\label{eq-irel-inclusion}
\alpha\subseteq\beta &\implies& R(\alpha)\subseteq
R(\beta) \\[1ex]
\label{eq-irel-sub}
R(\alpha\setminus\beta) &=& R(\alpha)\setminus
R(\beta)
\end{eqnarray}
\end{corollary}
\begin{proof}
\eqref{eq-irel-inclusion} is trivial (for any binary relation), and \eqref{eq-irel-sub} is immediate
from the properties above:
\begin{equation}
R(\alpha)\;\;=\;\;R\big((\alpha\setminus\beta)\uplus
(\alpha\cap\beta)\big)\;\;
\stackrel{(\ref{eq-stab-cap},\ref{eq-stab-uplus})}=\;\; R(\alpha\setminus\beta)\uplus \big(R(\alpha)\cap R(\beta)\big)
\end{equation}
hence
\begin{equation*}
R(\alpha\setminus\beta)\;\;=\;\; R(\alpha) \setminus \big(R(\alpha)\cap R(\beta)\big)
\;\;=\;\; R(\alpha)\setminus R(\beta)
\end{equation*}

\vspace*{-3.3ex}

\end{proof}

Refer once again to the diagram \eqref{diag-pullback}.
Recall that we identify a
a pair $\langle \alpha,\beta\rangle$ of subsets $\alpha\subseteq A$
and $\beta\subseteq B$ with $\alpha+\beta\subseteq A+B$.
Intersection, union and inclusion of synchronisations are defined via
this identification.
Write $h(\langle\alpha,\beta\rangle)=h(\alpha)$ and
$g(\langle\alpha,\beta\rangle)=g(\beta)$.
Thus $\sigma\subseteq A+B$ is a synchronisation iff
\begin{equation}\label{eq-hsigma-gsigma}h(\sigma)\;\;=\;\;g(\sigma)\end{equation}
\begin{lemma}\label{lemma-sync-closure}
Synchronisations are closed under union, intersection and
subtraction:
\begin{itemize}
\myitem{a} if\/ $S$ is a set of synchronisations then $\bigcap S$ and $\bigcup S$
are synchronisations;
\myitem{b} if $\sigma$ and $\tau$ are synchronisations then
$\sigma\setminus\tau$ is a synchronisation.
\end{itemize}
\end{lemma}
\begin{proof}
\begin{equation*}
g\left(\bigcap S\right)
\;\;\;\stackrel{\eqref{eq-stab-cap}}=\;\;\;
\bigcap_{\sigma\in S\!\!}g(\sigma)
\;\;\;\stackrel{\eqref{eq-hsigma-gsigma}}=\;\;\;
\bigcap_{\sigma\in S\!\!}h(\sigma)
\;\;\;\stackrel{\eqref{eq-stab-cap}}=\;\;\;
h\left(\bigcap S\right)
\end{equation*}
The $\bigcup$ and subtraction cases are analogous, via \eqref{eq-stab-cup} and \eqref{eq-irel-sub}.
\end{proof}
\begin{lemma}
Distinct paths are disjoint: if $\gamma,\gamma'\in P$ then
\begin{equation}\label{eq-paths-disjoint}
\gamma\neq\gamma'\;\;\implies\;\;\gamma\cap\gamma'=\emptyset
\end{equation}
\end{lemma}
\begin{proof}
$\gamma\cap\gamma'$ is a synchronisation by intersection-closure (Lemma~\ref{lemma-sync-closure}).
If $\gamma\neq\gamma'$ and $\gamma\cap\gamma'\neq\emptyset$ then
$\gamma\cap\gamma'$ is a synchronisation strictly smaller than at
least one of $\gamma$ or $\gamma'$, contradicting minimality.
\end{proof}
\begin{lemma}[Decomposition]\label{lemma-path-decomp}
Every synchronisation $\sigma$ is the disjoint union of its paths:
\begin{equation}\label{eq-path-decomp}
\sigma\;\;\;\;=\;\;\;\;\biguplus\;\big\{\,\gamma\subseteq\sigma\,:\,\gamma\text{ is a path}\:\big\}
\end{equation}
\end{lemma}
\begin{proof}
Paths are disjoint by the previous lemma, so it remains to show that
every link $c\in\sigma$ is in some (necessarily unique) path $\gamma_c$.
(Automatically $\gamma_c\subseteq\sigma$, by minimality with respect
to $\gamma_c\cap\sigma$.)
Define
\begin{equation}
\gamma_c\;\;\;=\;\;\;\bigcap\:\big\{\,\tau\,:\,\tau\text{ is a synchronisation and }c\in\tau\,\big\}\;,
\end{equation}
a synchronisation by intersection-closure
(Lemma~\ref{lemma-sync-closure}) and non-empty since it contains
$\sigma$.  We must show that $\gamma_c$ is minimal among all
non-empty synchronisations (not merely among those containing $c$).
Suppose $\mu\subsetneq\tau$ is a non-empty synchronisation.
Let $\overline\mu\,=\,\tau\setminus\mu$, a synchronisation by
subtraction-closure (Lemma~\ref{lemma-sync-closure}).
Then one of $\mu$ and $\overline\mu$ is a synchronisation containing
$c$ which is strictly smaller than $\gamma_c$, a contradiction.
\end{proof}
\begin{proofof}{Theorem~\ref{thm-pullback}}
The square \eqref{diag-pullback} commutes:
$$g(p\langle\alpha,\beta\rangle)\;\;\stackrel{\eqref{eq-p}}=\;\;g(\alpha)\;\;\stackrel{\eqref{eq-sync}}=\;\;h(\beta)\;\;\stackrel{\eqref{eq-q}}=\;\;h(q\langle\alpha,\beta\rangle)\:.$$
Suppose $\inlinelinking{A}{p'}{P'}{q'}{B}$ yields an analogous
commuting square: $gp'=hq'$.
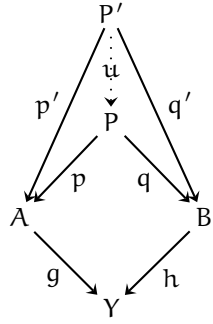
$$\newcommand{\hozgap}{\hspace{13ex}}\psset{arrows=->,labelsep=2pt}
\begin{array}{c}
\Rnode{P'}{P'} \\ [6ex]
\Rnode P P \\[5ex]
\Rnode A A \hozgap \Rnode{B}{B} \\[4.7ex]
\Rnode Y Y
\end{array}
\ncline{A}{Y}\nbput{g}
\ncline{B}{Y}\naput{h}
\ncline{P}{A}\naput{p}
\ncline{P}{B}\nbput{q} 
\ncline{P'}{A}\nbput{p'\!}
\ncline{P'}{B}\naput{q'}
\psset{linestyle=dotted}
\ncline{P'}{P}\ncput[npos=.54]{u}
$$
For $d\in P'$ let
\begin{equation}\label{eq-sigma-def}
\sigma(d)\;\;\;=\;\;\;p'(d)+q'(d)\;\;\;\subseteq\;\;\; A+B
\end{equation} 
which is a synchronisation since $gp'=hq'$.  Define $u:P'\to P$ by
taking $u(d)$ as the set of all paths within $\sigma(d)$:
\begin{equation}\label{eq-u-def}
u(d)\;\;\;=\;\;\;\{\,\gamma\in P\,:\,\gamma\subseteq \sigma(d)\,\}
\end{equation}

\emph{Claim: $u$ is injective.}
If $u(d)\cap u(e)\neq\emptyset$ there exists a path
$\gamma$
such that $\gamma\subseteq\sigma(d)\cap\sigma(e)$, say
$\gamma=\alpha+\beta$ with $\alpha\subseteq A$ and $\beta\subseteq B$.
Hence $\alpha\subseteq p'(d)\cap p'(e)$ and $\beta\subseteq q'(d)\cap
q'(e)$.
Since $\gamma$ is a path, it is non-empty, so $\alpha$ or $\beta$ is
non-empty, say $\alpha$.  Thus $d=e$ by injectivity of $p'$. \qed

\vspace{.8ex}\emph{Claim: $pu=p'$ and $qu=q'$.}
Suppose $a\in p'(d)$.
Let $\gamma$ be the unique path such that $a\in\gamma$ and
$\gamma\subseteq \sigma(d)$, existing by
Lemma~\ref{lemma-path-decomp}.  
Then $\gamma\in u(d)$ (by \eqref{eq-u-def}) and $a\in p(\gamma)$
(since $a\in A$ and $p$ projects subsets of $A+B$ to subsets of $A$),
hence $a\in p(u(d))$, so $p'\subseteq pu$.

Conversely, suppose $a\in p(u(d))$, \ie, there exists $\gamma\in P$
such that $a\in p(\gamma)$ and $\gamma\in u(d)$.  By \eqref{eq-u-def}
we have $\gamma\subseteq \sigma(d)$, so $a\in p(\sigma(d))$, by
\eqref{eq-irel-inclusion}.
Since $p(\sigma(d))=p'(d)$ (because $p$ projects)
we have $a\in p'(d)$.  Hence $pu\subseteq p'$.

Since $p'\subseteq pu$ and $pu\subseteq p$, we have $p'=pu$, whence
$q'=qu$, by symmetry.   \qed\vspace{.8ex}

Finally, we must prove that $u$ is unique, \ie, the commuting
triangles $pu=p'$ and $qu=q'$ determine $u$.
Let $\hat u:P'\to P$.
Given $r:A\to M$ and $s:A\to N$ write $[r,s]$ for the corresponding
injective relation $A\to M+N$.
Thus $p'=p\hat u$ and $q'=q\hat u$ iff $[p,q]\hat u = [p',q']\hat u$.
Paths are non-empty, so $[p,q]$ is total, hence monic
(Lemma~\ref{lemma-monic}).
Therefore $[p,q]\hat u=[p,q]u$ implies $\hat u=u$.
\end{proofof}

\small
\bibliographystyle{myalphaams}

\providecommand{\bysame}{\leavevmode\hbox to3em{\hrulefill}\thinspace}

\end{document}